\documentclass[12pt]{article}
 \usepackage{amsmath,amsfonts,theorem}
\input epsf
 \newcommand{\ntag}[1]{} 
\numberwithin{equation}{section}
 
 \newtheorem{prop}{Proposition}

 \theorembodyfont{\rmfamily}
 \newtheorem{dfn}{Definition}

 \newcommand{\qed}{\ifhmode\unskip\nobreak\fi\quad\ensuremath\square}

 \newcommand{\ga}{\gamma}
 
 \newcommand{\om}{\omega}

 \newcommand{\Si}{\Sigma}

 \newcommand{\PP}{\mathbb P}
 \newcommand{\C}{\mathbb C}

 \newcommand{\R}{\mathbb R}
 \newcommand{\Z}{\mathbb Z}

\begin{document}

 \title{Homological orthogonality of "symplectic"
 and "lagrangian" - corrected version}

  \author{Nik. A. Tyurin\footnote{BLTP JINR (Dubna), e-mail: ntyurin@theor.jinr.ru}}

\date{}

\maketitle

\begin{abstract} In this remark we discuss a relationship between
(co)homology classes of a symplectic manifold realized by
symplectic and lagrangian objects. We establish some transversality
condition for the classes, realized by symplectic divisors and smooth
lagrangian tori with some special condition on their
intersections.
 \end{abstract}
 \bigskip

Let $(M, \om)$ is a compact symplectic manifold equipped with symplectic form $\om$.
Then for the (co) homology classes of $M$ one can study the realizability problem;
the fixed symplectic structure distinguishes from the space of smooth realization
of homology classes two types of subobjects:

--- symplectic submanifolds

and

--- lagrangian submanifolds.

A submanifold $D \subset M$ is symplectic if the restriction $\om|_D$ is a symplectic form
on $D$; a submanifold $S \subset M$ is isotropic if the restriction $\om|_S$ is trivial and
$S$ is lagrangian if $\dim S = \frac{1}{2} \dim M$ thus lagrangian means maximal isotropic.

The situation studied in this remark is the following: we take a cohomology class $c
\in H^2(M, \Z)$ and a homology class $\Si \in H_n(M, \Z)$ such that the class $P.D.(c)
\in H_{2n-2}(M, \Z)$ is realized
by a smooth symplectic divisor $D \subset M$ and $\Si$ is realized by a smooth orientable
 lagrangian submanifold $S \subset M$. In what follows $\dim M = 2n$ therefore
$\dim D = 2n-2$ and $\dim S = n$.

Standard intersection theory says that two generic submanifolds of dimensions
$2n-2$ and $n$ in an ambient manifold of dimension $2n$ have the intersection
of dimension
$$
2n-2 + n - 2n = n-2,
$$
and it is the case of the {\it transversal} intersection. At the same time
for symplectic divisors and lagrangian submanifolds we have the following simple fact:
for any $D$ and $S$ at each point $p \in D \cap S$ of their intersection the dimension
of $T_p (D \cap S)$ is either $n-2$ or $n-1$ (of course, if $D \cap S$
is non empty). The proof is obvious --- the intersection $D \cap S$ must be isotropic in
$(D, \om|_D)$ therefore the dimension must be less or equal to $n-1$; the transversality
arguments show that it must be greater or equal to $n-2$.

If for given $D$ and $S$ the intersection $D \cap S$ has pure dimension $n-2$ this situation
is called transversal; thus due to the remark above it's quite natural to give a name
for the opposite situation:

\begin{dfn} Let $D \subset M$ is a symplectic divisor and $S \subset M$ is a lagrangian
submanifold of a symplectic manifold $M, \om$. Then we say that they have co-transversal
intersection if either at each point $p \in D \cap S$ the intersection has dimension $n-1$ or
$D \cap S$ is empty.
\end{dfn}

{\bf Example} Let $M$ be $2n$ - dimensional torus $T^{2n}$ decomposed into two tori
$T^2 \times T^{2n-2}$ with the product symplectic structure coming from symplectic structures
over $T^2$ and $T^{2n-2}$. Then if we take $D = p \times T^{2n-2}$ and $S = \ga \times S_1$
where $p$ is a point in $T^2$, $\ga$ is a circle in $T^2$ and $S_1$ is a lagrangian submanifold
in $T^{2n-2}$ then $D$ and $S$ has co- transversal intersection in $T^{2n}$.

{\bf Another example} In \cite{Aur} one studies a situation when some lagrangian
submanifold $S \subset M$ has 1 - dimensional projection to $\C \PP^1$ with respect
to a Lefschetz pencil
$$
f: M \to \C \PP^1.
$$
 This means that for almost all fibers of this Lefschetz
pencil the intersection with $S$ is co-tansversal. Since the strong result of \cite{Aur}
ensures that in dimension 4 for integral symplectic manifold $M, \om$ (this means that
$[\om] \in H^2(M, \Z)$) for any lagrangian submanifold one can find a Lefschetz pencil
with the desired property, the case of co-transversal intersection is not too rare.

The last example leads to a natural extension of the definition given above:

\begin{dfn} Let $f: M \to \C \PP^1$ is a Lefschetz pencil on a symplectic manifold $M$.
Then $f$ is co-transversal to a smooth  lagrangian submanifold $S \subset M$ if
each smooth fiber of $f$ has co-transversal intersection with $S$.
\end{dfn}

A model example can be found in \cite{Stef}: there one takes for a given
lagrangian submanifold $S \subset X$ of a projective algebraic manifold $X$
(with a fixed  Kahler metric of the Hodge type) some small deformation
$S_1 \subset X$ which is real analytic and then for a given Morse function
$F$ on $S_1$ one constracts an analytic extension $f$ which can be extended
to a Lefschetz pencil on the whole $X$. The crucial point is that
$f$ by the construction is co-transversal to $S_1$.

The main background problem seems to be interesting is the following: for a given
lagrangian submanifold $S \subset M$ of a given symplectic manifold  construct
a Lefschetz pencil of a given topological type which is co-transversal to
$S$. Equivalently,  construct a Lefschetz pencil $f: M \to \C \PP^1$ such that
$f(S)$ is a loop in $\C \PP^1$. Indeed, by the definition of the Lefschetz pencil,
see \cite{Don}, the map $f$ is defined on the complement of a symplectic
submanifold $N \subset M$ of real codimension 4. This means that the intersection
$N \cap S$ has real codimension at least 2 and in the co-transversal case
$f(S)$ has no boundary in $\C \PP^1$.

The property of co-transversal intersections can be illustrated
by the following fact:

\begin{prop} Let $(M, \om)$ is a symplectic manifold with trivial
torsion in (co)homology. Let $c \in H^2(M, \Z)$ and $\Si
\in H_n(M, \Z)$ are (co) homology classes such that $P.D.(c)$ is realized by
smooth symplectic divisors and $\Si$ is realized by  smooth
lagrangian tori. Then the restriction $c|_S$ is trivial if
there are exist  representatives $D$ and $S$ --- symplectic and lagrangian
 respectively ---
with co-transversal intersection.
\end{prop}

Indeed, for the situation of the Proposition one can use a version
of the classical Darboux - Weinstein theorem namely there exists a neihborhood
$O_{\epsilon} (S) \subset M$ which is symplectomorphic to a neiborhood of the zerosection
in $T^*S$ endowed with the standard symplectic structure such that $D$, restricted to
$O_{\epsilon}(S)$, goes to a subbundle $W$ of $T*S|_{S_1}$ where $S_1 \subset S$ is
an orientable smooth $n-1$ - dimensional submanifold in $S$ (which is $D \cap S$). One could roughly
say that $D$ corresponds to $T^*S_1$ but since $T^*S_1$ doesn't lie in $T^*S$
we must speak about a subbundle of rank $n-1$. This subbundle defines a vector field
on $S_1$ up to scaling so we can fix a global smooth vector field $Y$ on $S$ which
is nontrivial on $S_1$ and  is annihilated by any section of $W$. To prove the Proposition it
 remains to take a smooth function $F \in C^{\infty}(S, \R)$ which is constant on $S_1$ and extend it
to $T^*S$ using $Y$ as a linear function along the fibers. The hamiltonian vector field $X_{F_Y}$
of the extended function $F_Y \in C^{\infty}(T^*S, \R)$  deforms $S$ to $S'$ which is
 a smooth lagrangian torus which represents the same homology class, $[S'] = [S]$, and which
has trivial intersection with $D$. This means that the line bundle $L$ with $c_1(L) = c =
P.D.([D])$ is restricted  to $S$ topologically trivial.

It is not hard to modify  Proposition 1 for the case of Lefschetz pencils:

\begin{prop} Let $(M, \om)$ is a symplectic manifold with a Lefschetz pencil
$f: M_1 \to \C \PP^1$ where $M_1 = M - N$, $N$ is symplectic of codimension 4.
If $S$ is a smooth lagrangian torus $S \subset M$ co-transversal to $f$ then
$c|_S$ is trivial, where $c \in H^2(M, \Z)$ is the topological type of the
pencil.
\end{prop}

Indeed, take a smooth fiber $f^{-1}(p) \subset M_1$ and restore a compact symplectic submanifold
$D = f^{-1}(p) \cup N$ whose topological type is the type of the pencil. Then
either $D$ is co-transversal to $S$ or $N$ has topologically trivial intersection with
$S$. The first case is directly adressed to Proposition 1 while in the second one
we need to deform $S$ to some $S'$ such that new $S'$ has co-trasnversal intersection with
$D$ and  again reduce the problem to  Proposition 1.

\end{document}